\documentclass[reprint]{revtex4-1}

\usepackage{latexsym}

\usepackage{amsmath}
\usepackage{amsthm}
\usepackage{amssymb}
\usepackage{amsfonts}
\usepackage{color}
\usepackage{tikz, lipsum, lmodern}
\usepackage{graphicx}
\usepackage[rightcaption]{sidecap}
\usepackage{dcolumn}
\usepackage{bm}

\newtheorem{prop}{Proposition}[section]

\newtheorem{rem}[prop]{Remark}

\newcommand{\ba}{\begin{array}}
\newcommand{\ea}{\end{array}}

\newcommand{\be}{\begin{equation}}
\newcommand{\ee}{\end{equation}}

\numberwithin{equation}{section}

\begin{document}

\title{Computational enhancement of discrete gradient method}
\author{Artur Kobus\\
Uniwersytet w Bia\l ymstoku, Wydzia{\l} Fizyki\\
ul.\ Cio\l kowskiego 1L,  15-245 Bia\l ystok, Poland}
\email{a.kobus@uwb.edu.pl}

\date{\today}

\begin{abstract}
We propose new numerical approach to non-conservative dynamical systems. Our method being of low order, enhances qualitative performance of standard discrete gradient algorithm, thank to new concept of a reservoir. Paper is of explanatory character, focusing on concrete non-linear physical systems. Superiority of our new method with respect to standard discrete gradient method is observed.
\begin{description}
\item[Keywords]
Geometric Numerical Integration; Discrete Gradient Method; Non-conservative Systems
\item[PACS numbers]
45.10.-b, 02.60.Cb, 02.70.-c, 02.70.Bf
\end{description}
\end{abstract}
              
\maketitle      

\section{Introduction}

Non-trivial task of engaging the proper numercal algorithm for solving approximately target differential equation could be handled in two ways: using standard, general purpose integrators of high order as various Runge-Kutta-type methods, multi-step methods, or collocation methods, to name just a few \cite{I,MS}, or geometric methods exactly preserving qualitative properties like time-inversion symmetry, symplecticity or constants of motion \cite{B,B2,BC}, eventually reaching high order through composition and other tricks (see e.g. \cite{C,CK,CR3}).

We propose new method designed for the usual IVP
\be
\label{ivp}
\dot{\pmb{x}}=\pmb{f}(\pmb{x}), \qquad \pmb{x} (t_0) = \pmb{x}_0
\ee
treated in a very non-usual way. Here we picked system of autonomous form, but we could as well choose non-autonomous one.

We generalize already quite general approach of discrete gradient method \cite{MQR,QC} as constructed for Hamiltonian, Poisson and Lyapunov type systems, beyond the original scope. Especially we emphasize the application to non-linear, non-conservative dynamical systems. Convenient and effective treatment of these is possible thank to introducing novel kind of conserved quantity, the so-called ``computational invariant of motion".

Before we proceed to new results, we think a quick reminder on the classical discrete gradient method (only for the Hamiltonian case!) should be useful \cite{MQR}:
\begin{enumerate}
\item We discretize the problem on the time-grid $t_i=t_0+ih_i$,where $t_0$ is a chosen constant and $h_i$ is (possibly variable) time step.
\item We use the scheme
\be
\ba{l}
\frac{x_{i+1}-x_i}{h_i} = \frac{H(x_{i+1},y_{i+1})-H(x_{i+1},y_i)}{y_{i+1}-y_i},\\
\frac{y_{i+1}-y_i}{h_i} = \frac{H(x_{i+1},y_i)-H(x_i,y_i)}{x_{i+1}-x_i}
\ea
\ee
to model discrete analogue of Hamiltonian equations.
\item Cross-multiplying above equations leads to instantaneous conclusion
\be
H(x_{i+1},y_{i+})=H(x_i,y_i).
\ee
\end{enumerate}

The paper is organized as follows: in sec. 2 we describe new approach and it technicalities, and also present new numerical algorithm. In Sec. 3 we analyze test-bench example of the damped harmonic oscillator, being the simplest dissipative system, that posses well-known and easily-available exact solution, hence we can easily estimate the order of errors and search for eventual advantages with respect to other low-order integrators. In sec. 4 we pass on to important examples of Duffing and Van der Pol oscillators, illustrating features of new approach. This helps us to derive number of conclusions and future plans, contained in sec. 5, with which the paper ends.

\section{Computationally conserved quantities and new numerical scheme}

Let us consider a basic example
\be
\label{eom}
\ba{l}
\dot{x}=y,\\
\dot{y}= -x - D(x,y),
\ea
\ee
which is essentially harmonic oscillator with some damping/forcing included, permitting non-conservative behavior to occur (dependence on $y$ is crucial-otherwise we still get Hamiltonian case). Pure harmonic oscillator would possess a constant of motion $H(x,y)=\frac{1}{2} x^2 + \frac{1}{2} y^2$, but due to (\ref{eom}) there is
\be
\dot{H} = - y D(x,y).
\ee
Now we search for a quantity with exactly opposite behavior of the time derivative, in other words we need a function $z$ such, that
\be
\dot{z} = y D(x,y)
\ee
and of course, we can formally solve above equation with the integral
\be
\label{w}
z= - \int D(x,y) y dt,
\ee
since $(x,y)$ are of course functions of independent variable $t$.

With definition of computationally conserved quantity
\be
K=H+z
\ee
we can directly check it is conserved.

Note that $z$ nor $K$ are not well-defined functions, since they depend on the path of integration in (\ref{w}). Moreover, to find $z$, we would need to explicitly know $x$ and $y$ as functions of time. Hence the anticipated redundancy of $z$ variable. We will call $z$ the \emph{reservoir variable}. 

New scheme is of the form:
\be
\label{egr}
\ba{l}
\frac{x_{i+1} -x_i}{h} = \frac{K(x_{i+1},y_{i+1},z_{i+1})-K(x_{i+1},y_{i},z_{i+1})}{y_{i+1}-y_i},\\
\frac{y_{i+1} -y_i}{h} = - \frac{K(x_{i+1},y_{i},z_{i+1})-K(x_{i},y_{i},z_{i})}{x_{i+1}-x_i},\\
w_{i+1} = w_i + D(x_i^*,y_i^*) (x_{i+1}-x_i),
\ea
\ee
the last equation appears since $y dt = dx$, due to (\ref{eom}). $(x^*,y^*)$ is a point calculated arbitrarily, but consequently throughout solving the equations of motion. Cross-multiplying first two equations of above scheme yields
\be
K(x_{i+1},y_{i+1},z_{i+1})=K(x_i,y_i,z_i).
\ee

\section{Example of the simple harmonic oscillator}

Having established new tools, we begin testing with simple example
\be
\ba{l}
\dot{x}=y,\\
\dot{y} = -x -b y,
\ea
\ee
where $b$ is the damping parameter.

As in former section, we adjoin the reservoir variable to the system
\be
z=b \int y^2 dt,
\ee
so that $K=H(x,y)+z$ is conserved due to equations of motion.

We turn the system into a discrete one
\be
\ba{l}
\frac{x_{i+1} -x_i}{h} = \frac{1}{2} (y_i+y_{i+1}),\\
\frac{y_{i+1} -y_i}{h} = - \frac{1}{2} (x_i+x_{i+1}) - b y_i^*,\\
z_{i+1} - z_i = b y_i^* (x_{i+1}-x_i)
\ea
\ee
where we will use $y_i^*=\frac{1}{2}(y_i+y_{i+1})$ for simplicity.

Additional quantity of which error we will investigate is the energy decrement
\be
R=\frac{H(x_{i+1},y_{i+1})}{H(x_i,y_i)}=\frac{K(x_0,y_0)-z_{i+1}}{K(x_0,y_0) - z_i},
\ee
this second form designed for en-GR scheme we use the advantage of having reservoir at our disposal. We accept initial conditions $x_0=1.3, y_0=-2.2$.

In this paper all numerical experiments are performed with simulation period $T=100.0$, the time-step $h=0.001$, and accuracy tolerance $\varepsilon = 10^{-16}$ while solving non-linear algebraic equations by iteration. We use explicit Euler scheme as a predictor.

We xompare en-GR, st-GR (standard discrete gradient), IMR and SV. As a substitute for exact solution (if it is not known) we use numerical sequences generated by $^{th}$ order Runge-Kutta algorithm ($3/8$ rule) with a way tinier time-step $h_{RK} = 10^{-6}$. Also we substitute midpoint approximation in place of $x_i^*,y_i^*$.

\begin{minipage}[c]{0.495\textwidth}

\includegraphics[width=6cm, height=5cm]{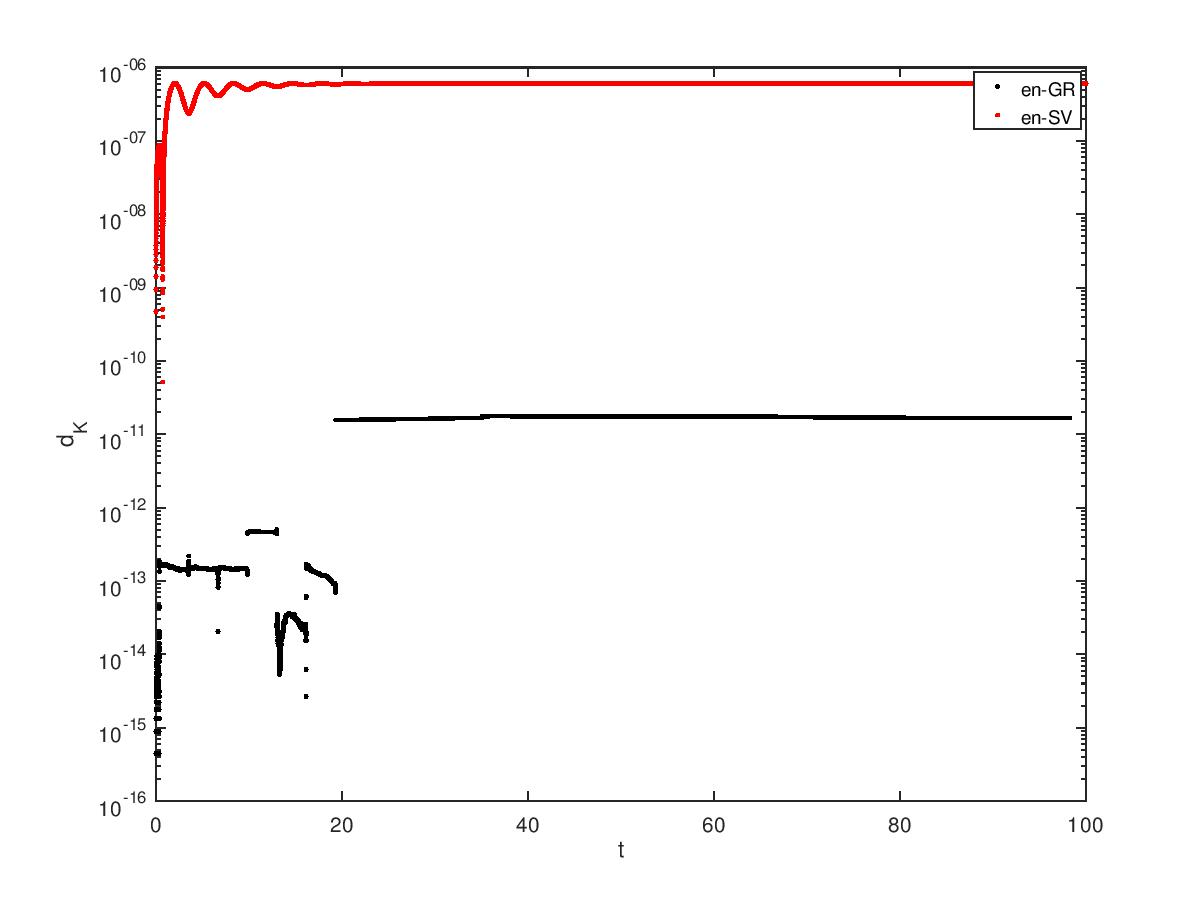}

\footnotesize{Figure 3.1.: Deviation from initial value of $K$.\newline We adjoined reservoir to SV to see, how it conserves 
$K$.}
$\\[1ex]$

\end{minipage} \begin{minipage}[c]{0.01\textwidth}   \end{minipage} \begin{minipage}[c]{0.495\textwidth}

\includegraphics[width=6cm, height=5cm]{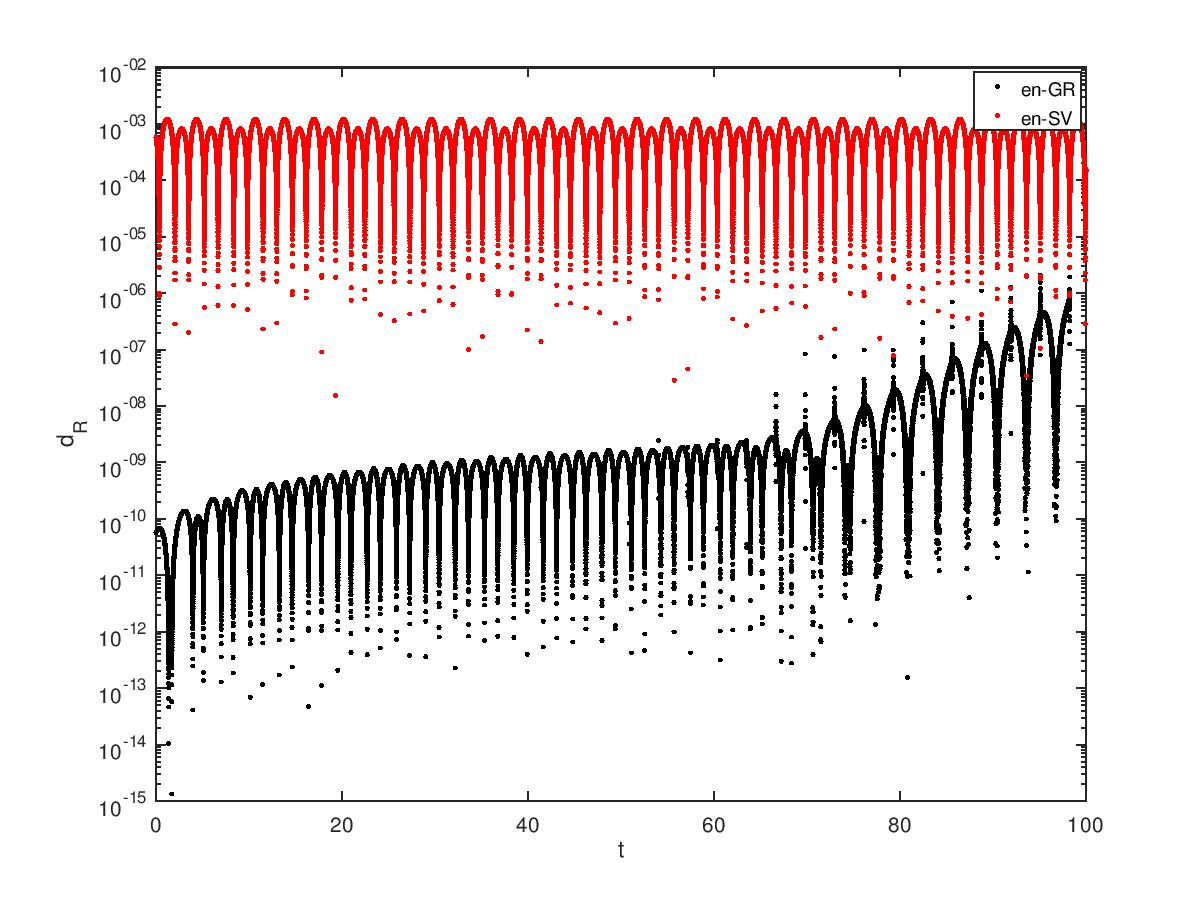}

\footnotesize{Figure 3.2.: Comparison of deviation from theoretical value of energy decrement. en-GR is better through six orders of magnitude. IMR yields the same results as SV.  $b=0.2$.}
$\\[1ex]$

\end{minipage}

\section{Application to non-linear systems}

\subsection{Van der Pol oscillator}

Equations of motion of Van der Pol oscillator are
\be
\ba{l}
\dot{x}=y,\\
\dot{y} = - x +a(1-x^2)y,
\ea
\ee
what in discrete version with a reservoir, reads
\be
\ba{l}
x_{i+1} = x_i +\frac{h}{2} (y_i+y_{i+1}),\\
y_{i+1} = y_i - \frac{h}{2} (x_i+x_{i+1} - 2a (1-(x_i^*)^2)y_i^*),\\
z_{i+1} = z_i - a(1-(x_i^*)^2) y_i^* (x_{i+1}-x_i)
\ea
\ee
where we picked $a=1.0$.

With midpoint approximation in place of $x_i^*,y_i^*$, en-GR becomes the same algorithm as IMR, but with computational conserved quantity $K$. Standard gradient algorithm cannot be used in this case. We apply initial conditions $x_0=3.42,y_0=2.5$.

\subsection{Damped Duffing oscillator}
\label{duf}

The continuous system is
\be
\ba{l}
\dot{x}=y,\\
\dot{y} = x - x^3 -b y,
\ea
\ee
so for the conservative case Hamiltonian reads
\be
H(x,y) = \frac{1}{2} y^2- \frac{1}{2} x^2 + \frac{1}{4} x^4.
\ee
With addition of a reservoir the enhanced discrete gradient scheme (\ref{egr}) reads
\be
\ba{l}
x_{i+1}=x_i + \frac{h}{2} (y_i+y_{i+1}),\\
y_{i+1}=y_i +\frac{h}{2} ((x_i+x_{i+1}) (1-\frac{x_i^2+x_{i+1}^2}{2}) + 2 b y_i^*),\\
z_{i+1} = z_i+ b y_i^* (x_{i+1}-x_i)
\ea
\ee
with $b=0.2$ again. This time $x_0= -6.0, y_0= 2.5$.

\begin{rem}
\label{en-bet}
Remembering that we substitute simple midpoint approximation for $x_i^*,y_i^*$, we achieve the same form of en-GR and the standard gradient algorithm (st-GR).
\end{rem}

All tested schemes (with exception of st-GR) point towards left minimum and this is the reason why $x$-error of this methodn is so huge when compared to en-GR, however we see trouble begins much earlier. Reminding now remark (\ref{en-bet}) we have natural and clear understanding that gradient scheme not equipped with a conserved quantity can loose proper pace of energy decreasing.

\begin{minipage}[c]{0.495\textwidth}

\includegraphics[width=6cm, height=5cm]{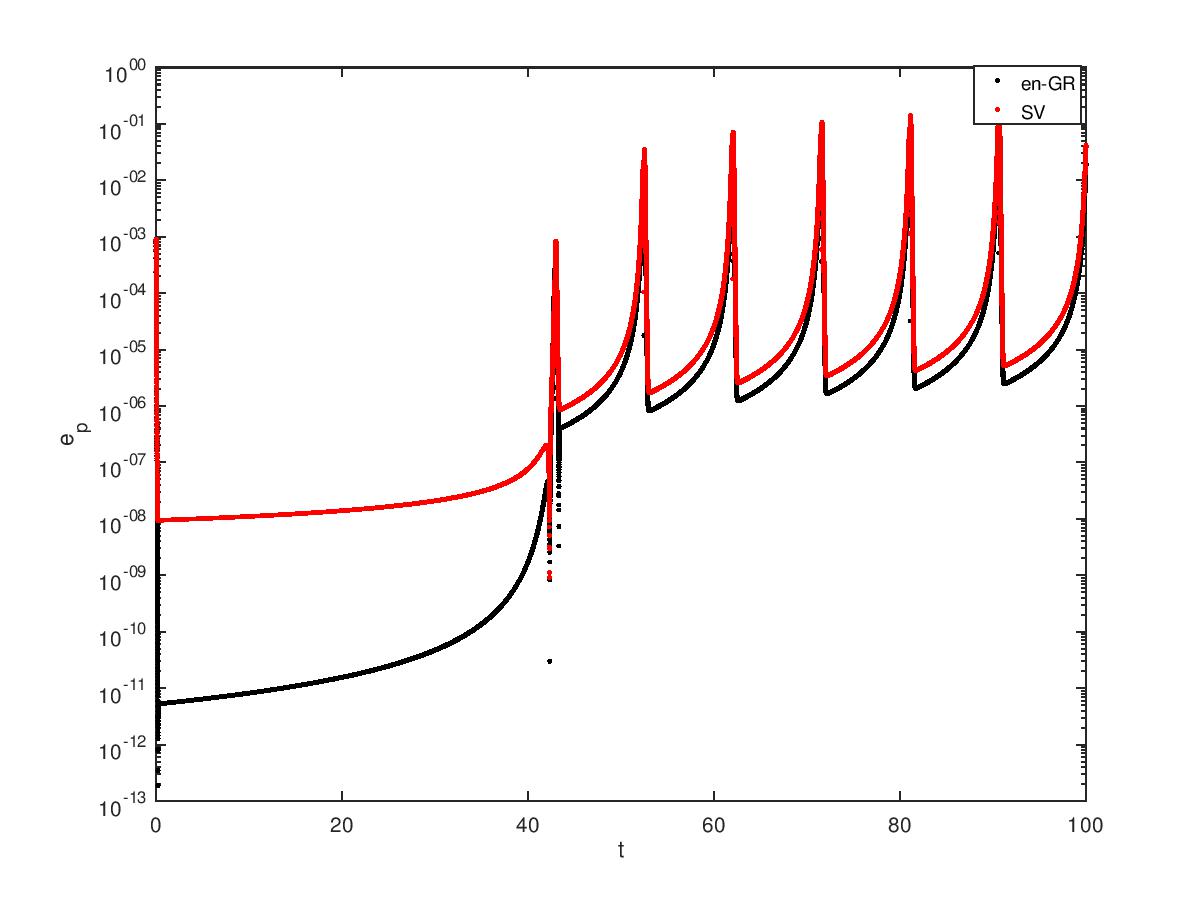}

\footnotesize{Figure 4.1.: Error of $p$ variable as measured with respect to RK4 solution in the case of Van der Pol oscillator.}
$\\[1ex]$

\end{minipage} \begin{minipage}[c]{0.01\textwidth}   \end{minipage} \begin{minipage}[c]{0.495\textwidth}

\includegraphics[width=6cm, height=5cm]{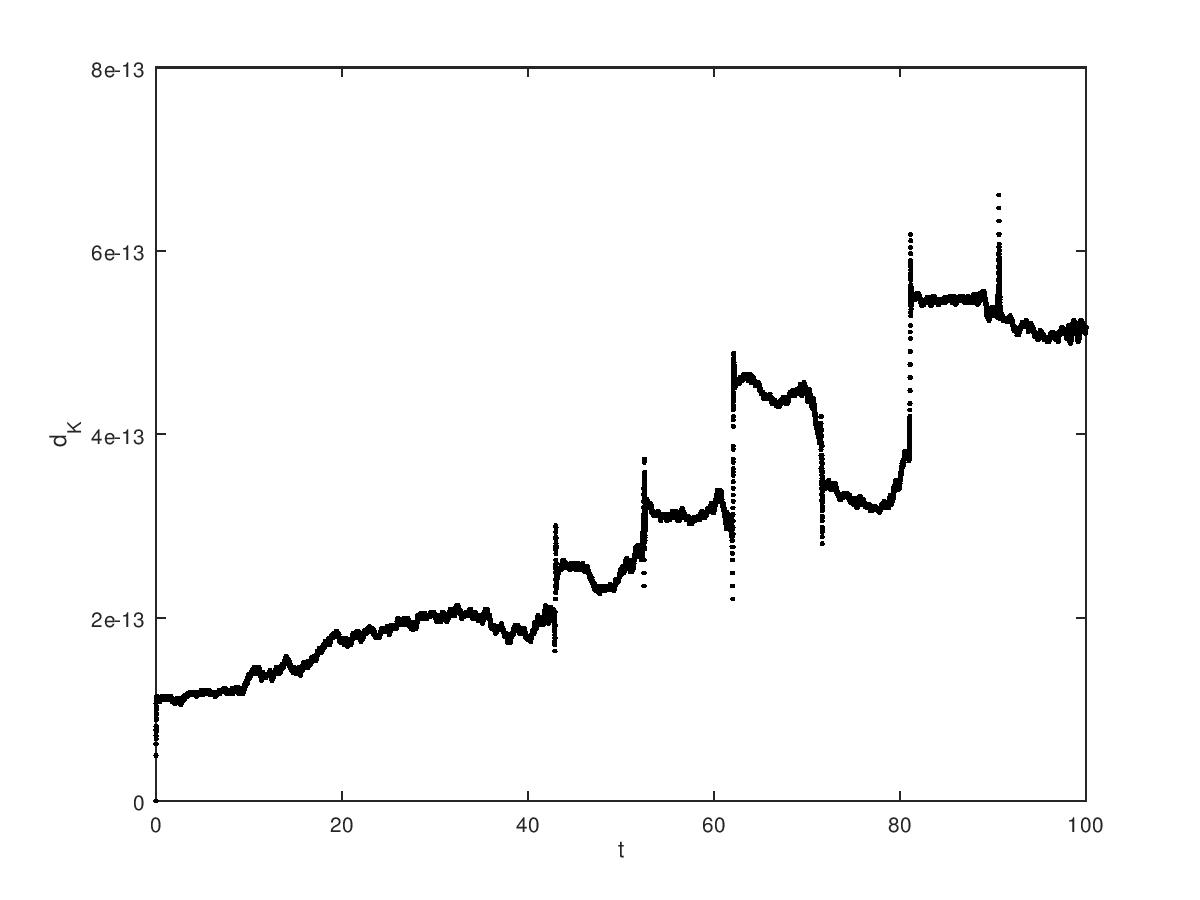}

\footnotesize{Figure 4.2.: Deviation from value of $K$ for en-GR (Van der Pol).}
$\\[1ex]$

\end{minipage}

\begin{minipage}[c]{0.495\textwidth}

\includegraphics[width=6cm, height=5cm]{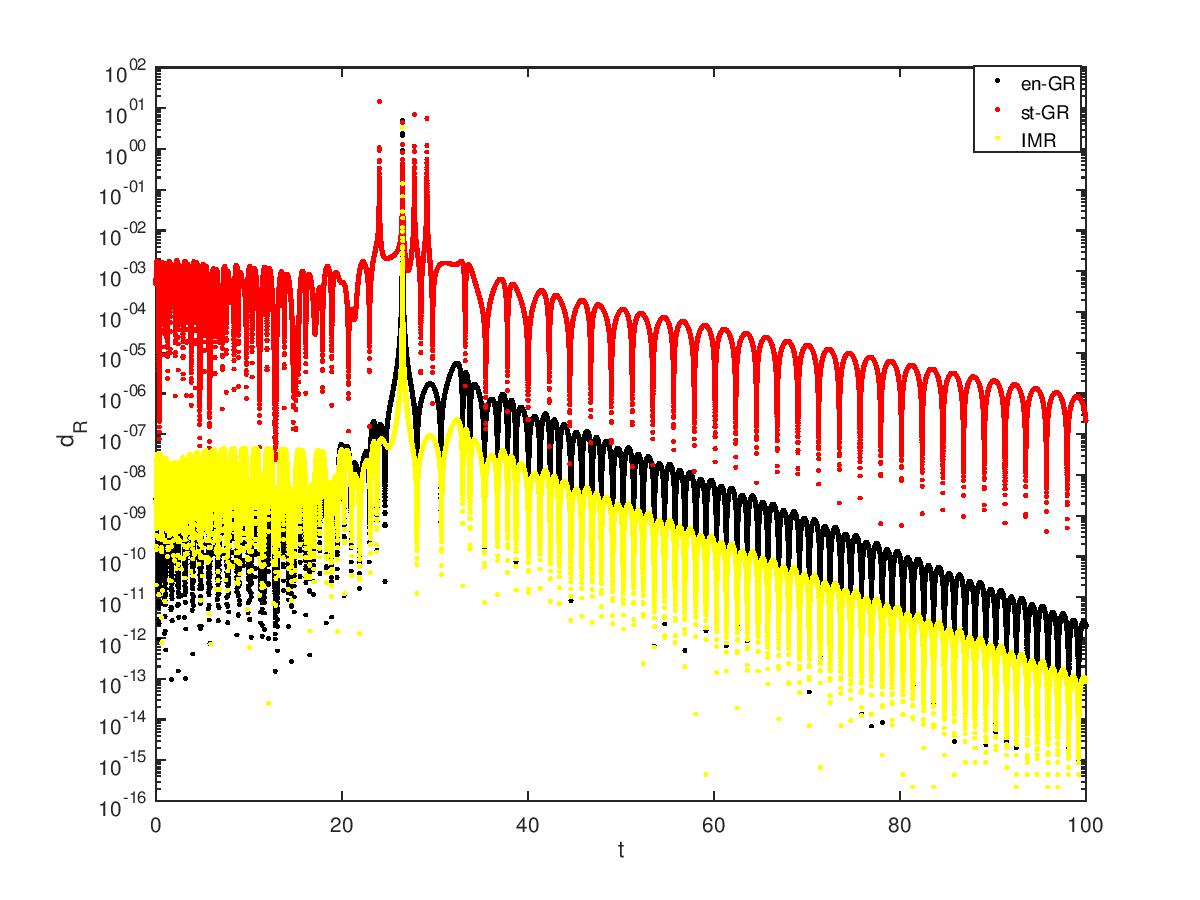}

\footnotesize{Figure 4.3.: Deviation from the energy decrement.\newline SV performing equally with IMR, hence it is omitted (Duffing).}
$\\[1ex]$

\end{minipage} \begin{minipage}[c]{0.01\textwidth}   \end{minipage} \begin{minipage}[c]{0.495\textwidth}

\includegraphics[width=6cm, height=5cm]{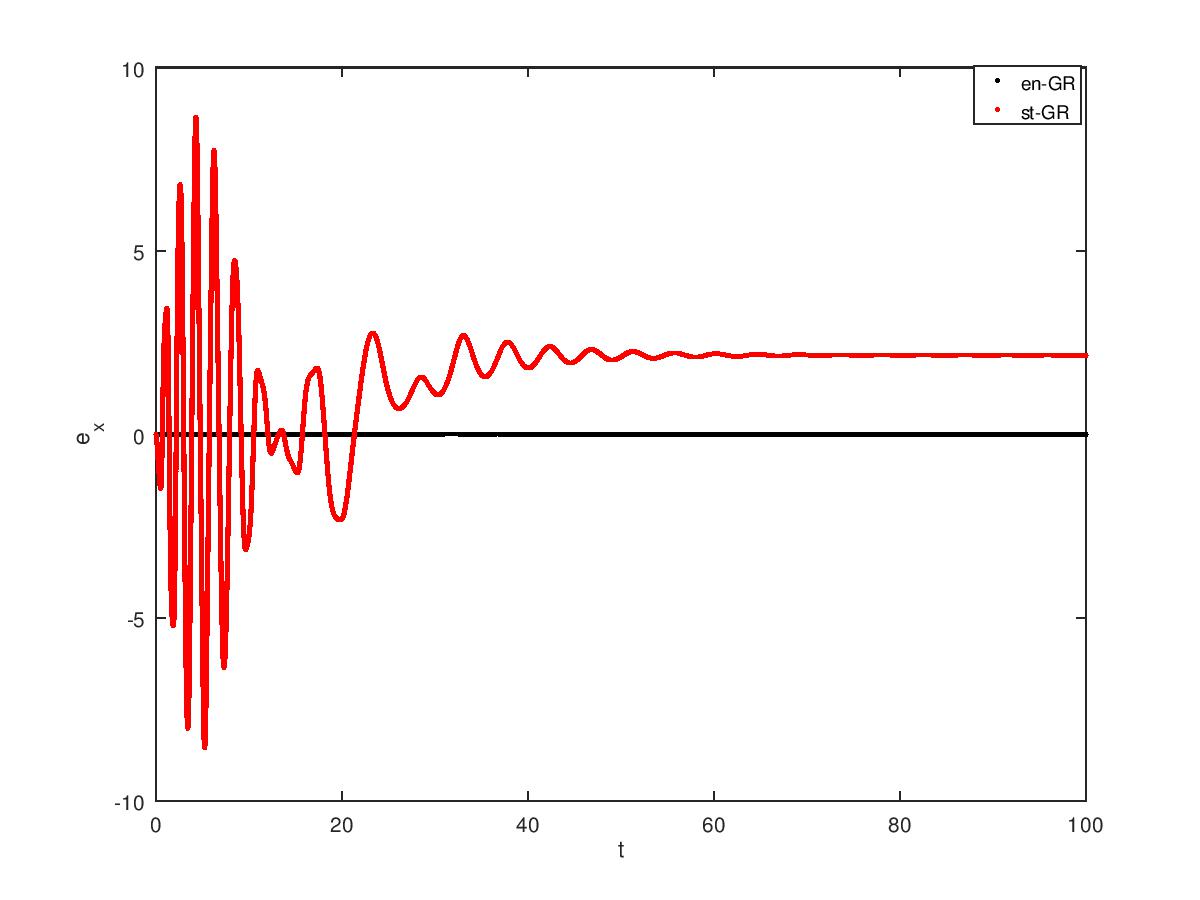}
\footnotesize{\newline Figure 4.4.:  Global error of $x$, st-GR fails\newline qualitatively, landing in the right minimum (Duffing).}
$\\[1ex]$

\end{minipage}

\section{Conclusions and plans}

We introduced new numerical schemes based on novel concepts of reservoir and computationally conserved quantity. This new algorithms belong to geometric numerical integration methods, preserving some quantities by in-built mechanisms. 

We performed multiple numerical tests to check if our approach yields results compatible with other well-known integrators. Enhanced discrete gradient scheme presents second order behavior and is especially useful when the energy decrement plays   important role. Morever, it enables us to use discrete gradient procedure even if there is no constant of motion nor Lyapunov function for the system.

Future perspective of research in this direction is threefold:
\begin{enumerate}
\item Main focus of this paper was to briefly introduce new methods with some illustrative examples; we skipped theoretical analysis of the scheme. In fact, basic calculations suggest our approach is of close oeigin to \cite{MM},
\item As shown in section \ref{duf}, enhanced discrete gradient has shown some advantages in preserving trajectories of dynamical systems. There were performed numerical simulations, that further confirm this property of en-GR (good example gives the heteroclinic motions of dynamical systems). Moreover, we did not search for the optimal form of $x_i^*,y_i^*$, also interesting would be to test long-time performance of this method,
\item Higher dimensional counterpart of this method is also in preparation. It seems like matrix formulation of \cite{MQR} will  apply with minor changes in this case.
\end{enumerate}

\end{document}